\documentclass{elsart}
\usepackage{natbib,amsmath,amssymb,latexsym,url,graphicx}
\usepackage{graphics}
\usepackage{psfig}
\usepackage{color}
\def\para{\vspace{3mm}}

\def\deg{{\rm deg}}

\def\qed{\hfill  \framebox(5,5){}}
\def\lcoeff{{\rm lcoeff}}

\def\Res{{\rm Res}}

\newtheorem{theorem}{{\bf Theorem}}
\newtheorem{remark}{{\bf Remark}}
\newtheorem{definition}[theorem]{{\bf Definition}}

\newtheorem{proposition}[theorem]{{\bf Proposition}}

\newtheorem{example}{{\bf Example}}

\begin{document}
\begin{frontmatter}



\title{Applications of Level Curves to Some Problems on Algebraic Surfaces}

\author[a]{Juan Gerardo Alcazar\thanksref{proy}},
\ead{juange.alcazar@uah.es}
\author[a]{Juan Rafael Sendra\thanksref{proy}}
\ead{rafael.sendra@uah.es}

\address[a]{Departamento de Matem\'aticas, Universidad de Alcal\'a,
E-28871-Madrid, Spain}

\thanks[proy]{Both authors supported by the Spanish `` Ministerio de
Educaci\'on y Ciencia" under the Project MTM2005-08690-C02-01.}

\begin{abstract}
In \cite{JGRS}, a result to algorithmically compute the topology
types of the level curves of an algebraic surface, is given. From
this result, here we derive applications based on level curves to
determine some topological features of surfaces (reality,
compactness, connectivity) and to the problem of plotting.
\end{abstract}
\end{frontmatter}

\section{Introduction}

The study of the level curves of an implicit algebraic surface (i.e.
the sections of the surface with real planes parallel to the
$xy$-plane) gives a clue on how the surface is like. Take for
example the well-known example of the Whitney Umbrella, whose
equation is $x^2-y^2z=0$. It is clear that for $z>0$ the level
curves consist of two intersecting lines; for $z=0$, the level curve
reduces to one line; and for $z<0$, the level curves consist of just
one real isolated point. From this simple analysis one may get a
good mental picture of the surface. Moreover, it provides an initial
analysis to the topology of the surface, which can be used later to
provide topologically-correct plottings. For instance, in the case
of the Whitney Umbrella one has that a topologically-correct
plotting must show the so-called ``handle" of the surface, i.e. the
half-line corresponding to $x=0,y=0,z\leq 0$. Moreover, in addition
to the shape of the level curves, one also has the $z$-intervals
corresponding to each topology type.

The problem of computing the different topology types arising in the
family of level curves of a given algebraic surface (together with
the $z$-values corresponding to each type) has been addressed, for
instance in \cite{JGRS} and \cite{Mourrain-2}, where algorithms to
compute these topology types are given. Moreover, in \cite{JGRS}
applications of this approach to the study of offset curves  are
also given; see also \cite{JG-Sendra07} for further research on this
topic. Now, in this paper, we apply these results  on level curves
to determine topological information on the surface, which may be of
help for computing reliable plottings from the topological point of
view.

In this sense, the first application that we consider is to
algorithmically decide whether a given implicit algebraic surface
is real. That is, we provide an algorithm to check whether the
intersection of the surface with ${\mathbb R}^3$  is a
two-dimensional set (in the Euclidean topology) or not.
Furthermore, in the negative case, the algorithm can also be used
to analyze whether the real part of the surface is empty (like for
example $x^2+y^2+z^2+1=0$), or it consists of finitely many points
(the case of $x^2+y^2+z^2=0$, whose real part is the origin) or it
is a space curve (the case of $x^2+y^2=0$). This kind of surfaces
whose real part is not 2-dimensional may get error messages when
one tries to draw a plotting.

The second application concerns to the compactness of the surface.
In this case, since one works with implicit algebraic surfaces,
which are closed over the usual Euclidean topology, the algorithm
essentially checks whether the surface is bounded w.r.t. the
variables $x,y,z$, respectively.

The third application directly concerns to the problem of plotting.
In order to draw a plotting of a surface, the user has to introduce
as an input a ``box" $I=[a_1,a_2]\times [b_1,b_2]\times
[c_1,c_2]\subset {\Bbb R}^3$, so the output shows the part of the
surface lying inside the box. Now, if the user is interested in
computing a plotting where the main topological features of the
surface are shown (i.e. which makes clear  how the surface is like),
some previous information must be known in order to properly choose
$I$. Using the information on the topology of the level curves of
the surface w.r.t. the variables $x,y,z$, we provide an algorithm to
compute a box of this kind.

The fourth application is a symbolic-numeric algorithm to compute
the connected components of a surface by using level curves. Here,
the algorithm determines how the level curves are connected to each
other, and as a consequence also provides the number of connected
components of the surface. Since the algorithm in fact computes how
to join sections of the surface, we think that it might be used for
plotting purposes. This question might be explored in a future work.

The structure of the paper is the following. The first section
contains some preliminary notions and related results on level
curves. The second section is devoted to the problem of checking
whether a given implicit algebraic surface is real, or not. The
third section analyzes compactness. The problem of computing a
suitable box for plotting a surface (so the output shows the more
relevant topological features of the surface) is addressed in the
fourth section. The fifth section concerns to connectivity.

\section{Preliminaries on Level Curves}\label{Preliminaries}

In this section we provide some preliminary notions and results
concerning to level curves, that are taken from \cite{JGRS}. Thus,
we refer the interested reader to \cite{JGRS} for further
information. Moreover, here we also fix the notation to be used
along the paper, together with some hypotheses to be requested on
the surface $S$.

In the sequel, we consider an algebraic surface $S$ defined by an
square-free polynomial $F\in {\Bbb R}[x,y,z]$ with no univariate
factor only depending on the variable $z$, i.e. $S$ has no component
which is a plane parallel to the $xy$-plane. In this situation, the
{\sf level curves} of $S$ are the (plane) curves which are obtained
by intersecting $S$ with planes normal to the $z$-axis. Furthermore,
given $b\in {\Bbb R}$ we will denote the level curve corresponding
to the plane $z-b=0$ by $S_b$.

Note that one might similarly define the level curves corresponding
to the $x$-axis and the $y$-axis, respectively. Thus, when necessary
we will speak of $\xi$-level curves, where $\xi \in \{x,y,z\}$. Note
that the topology type of a $\xi$-level curve can be described by
computing a graph homeomorphic to it, which is a well-studied
problem (see for example  \cite{Lalo}, \cite{Hong} and many others).
Essentially, the vertices of this graph are the real critical points
of the curve (i.e. the real singular points and the real points
where the tangent is vertical), and the edges correspond to the
branches of the curve joining two vertices.

From Hardt's Semi-Algebraic Triviality Theorem (see \cite{Basu}) it
can be derived that the number of topology types of the level curves
of $S$, is always finite. In case that $S$ is compact and
non-singular the problem of determining these topology types can be
solved by using Morse Theory (see\cite{Basu}, \cite{Milnor}). In the
more general case of singular surfaces, two approaches can be
considered. The first one comes from Differential Topology and uses
elements of Whitney Stratification Theory (see \cite{Goresky}). This
approach has been used in \cite{Mourrain-2}. The second one comes
from Computer Algebra and uses as an essential tool the notion of
{\sf delineability} (see \cite{MacCallum}). This second approach has
been developed in \cite{JGRS}. In the rest of the section, we recall
some notions and results in \cite{JGRS} related to the computation,
by means of \cite{JGRS}, of the topology types of the level curves
of $S$.

\begin{definition} \label{Critical-Level}
We say that $a\in {\Bbb R}$ is a {\sf Critical Level Value} if the
topology of the level curves of $S$ changes at $z=a$, i.e.
$\forall\,\,$ $\epsilon>0$ there exists $a_{\epsilon}\in
(a-\epsilon,a+\epsilon)$ such that the level curves corresponding to
$a$ and $a_{\epsilon}$, have different topology types. Moreover, we
say that ${\mathcal A}\subset {\Bbb R}$ is a {\sf Critical Level
Set} of $S$, if ${\mathcal A}$ is finite and it contains all the
critical level values of $S$.
\end{definition}

Note that since the number of topology types of the level curves of
$S$ is always finite, the number of critical level values is also
finite and therefore a critical level set always exists. Moreover,
once  a critical level set ${\mathcal A}$ has been computed, the
topology types of the level curves can be obtained. Indeed, writing
${\mathcal A}=\{ \alpha_1,\ldots,\alpha_r\}$, we can decompose the
$z$-axis as
\[ (-\infty,\alpha_1)\cup \{\alpha_1\} \cup (\alpha_1,\alpha_2)\cup \cdots \cup
(\alpha_{r-1},\alpha_r)\cup \{\alpha_r\} \cup (\alpha_r,\infty).
\] Thus, taking a $z$-value for each open interval of the
partition, and applying the existing algorithms for computing the
topology type of a plane algebraic curve (see \cite{Lalo},
\cite{Hong}), one determines the topology type for all the
$z$-values corresponding to the interval. The remaining finitely
many level curves, corresponding to $F(x,y,\alpha_i)$,
$i=1,\ldots,r$, are also analyzed with the same strategy.

Therefore, the problem of determining the topology types of the
level curves of $S$ reduces to the computation of a critical level
set. Let us recall the results in \cite{JGRS} which allow to do
this. For this purpose, throughout this paper, besides of the
hypotheses required above (i.e. $S$ is implicitly defined by an
square-free polynomial $F\in {\Bbb R}[x,y,z]$ having no univariate
factor only depending on $z$), we {\bf also assume} that  the
leading coefficient of $F$ w.r.t.  $y$ does not depend on the
variable $x$. Note that this requirement can always be fulfilled by
applying if necessary a rotation around the $z$-axis, which does not
modify the topology of the level curves of $S$.

Now, we consider the following notation:  $D_w(G)$ denotes the
discriminant of a polynomial $G$ w.r.t. the variable $w$, i.e.
$D_w(G)=\Res_w(G,\frac{\partial G}{\partial w})$, $\sqrt{G}$ denotes
the square-free part of a polynomial $G$, and:
\[
\begin{array}{l}
M(x,z):=\left\{\begin{array}{cll} 0 & \,\, & \mbox{if }\deg_y(F)=0
\\ \sqrt{D_y(F)} & \,\, & \mbox{otherwise}
\end{array} \right.\\
R(z):=\left\{\begin{array}{cll} 0 & \,\, & \mbox{if }\deg_x(M)=0
\\ D_x(M(x,z)) & \,\, & \mbox{otherwise}
\end{array} \right.
\end{array}
\]

Then, the following result holds (see \cite{JGRS}):

\begin{theorem} \label{main-theorem}
Let $S$ satisfy the above hypotheses. Then, it holds that:
\begin{itemize}
\item [(1)] If $R(z)$ is not identically zero, then the set of
real roots of $R(z)$ is a critical level set.
\item[(2)] If $R(z)$ and $M$ are identically zero, then  the set of real
roots of $D_x(F)$ is a critical level set.
\item[(3)] If $R(z)$ is identically zero but $M$ is not
identically $0$, then the set of real roots of $M(z)$ is a critical
level set.
\end{itemize}
\end{theorem}

\begin{remark} \label{constant}
If $R(z)$ is a non-zero constant, then there is just one topology
type for all the level curves. Similarly for the case when $M(z)$ is
a non-zero constant.
\end{remark}



\section{First Application: Reality of Algebraic
Surfaces}\label{reality}

In this section we deal with the problem of deciding whether an
algebraic surface is real of not; i.e. whether $S\cap {\Bbb R}^3$
has dimension two  with the usual Euclidean topology. For this
purpose, we will show that the problem of checking the reality of
$S$ can be reduced to checking whether finitely many level curves of
$S$ are real, so the problem in ${\Bbb R}^3$ can be reduced to a
problem in ${\Bbb R}^2$. Observe that in order to check whether an
algebraic curve is real, one can use well-known algorithms, like for
example \cite{SW99}. The following theorem, which can be found in
\cite{Lang}, is essential. Here, we consider the usual definition of
regular point of a surface $F(x,y,z)=0$, i.e. a point $P$ of the
surface is {\it regular} if some first partial derivative of $F$ $P$
does not vanish at $P$.

\begin{theorem}\label{th-Lang}
$S$ is a real surface if and only if it has at least one real
regular point.
\end{theorem}

{\bf Proof.} See Theorem XI.3.6 in \cite{Lang}.

From this theorem, one may derive another result, concerning to the
level curves of $S$, which can be used to algorithmically check
whether $S$ is real. In order to see this, we need the following
previous result. We denote the partial derivatives of $F$ w.r.t. the
variables $x,y,z$, respectively, by $F_x$, $F_y$, $F_z$.

\begin{proposition} \label{prop-1} Let $F\in {\Bbb
R}[x,y,z]$, with  $\deg_y(F)>0$,  be the defining polynomial of the
surface $S$, and let ${\mathcal A}$ be the critical level set of $S$
computed by means of Theorem \ref{main-theorem}. If $a\in {\Bbb R}$
verifies some of the following conditions:
\begin{itemize}
\item[(i)] $a$ is a root of the leading coefficient of $F$ w.r.t.
$y$,
\item[(ii)] $a$ is a root of the leading coefficient of $M(x,z)$
w.r.t.  $x$,
\item[(iii)] the polynomial $F(x,y,a)$ has multiple factors,
\end{itemize}
then $a\in {\mathcal A}$.
\end{proposition}

{\bf Proof. }In order to prove the result, we distinguish the cases
when the polynomial $R$, defined in Section \ref{Preliminaries}, is
identically $0$, or not. Let us see first the case when $R\neq 0$.
Here, we observe that if $a$ is a root of the leading coefficient of
$F$ w.r.t. $y$, then by using the Sylvester form of the resultant
$\Res_y(F,F_y)$, one has that $z-a$ is a factor of
$D_y(F)=\Res_y(F,F_y)$ and therefore of $M$. Thus, it is also a
factor of the leading coefficient of $M$ w.r.t. $x$, and
consequently, (i) implies (ii). Now, let us see that the condition
(ii) implies that $R(a)=0$, and therefore $a\in {\mathcal A}$.
Indeed, using again the Sylvester form of the resultant
$\Res_x(M,M_x)=D_x(M)=R(z)$, one deduces that if $z-a$ is a factor
of the leading coefficient of $M$ w.r.t. $x$, then it is also a
factor of $R(z)$. Thus, $R(a)=0$ and $a\in {\mathcal A}$. Finally,
let us see that if $a$ verifies (iii) but not (i), then $R(a)=0$.
Let $H(x,z)=D_y(F)$; then $M(x,z)=\sqrt{H(x,z)}$. Since $a$ does not
verify (i),  the discriminant of $F$ w.r.t. $y$ behaves properly
under specializations when $z=a$, i.e. $D_y(F_a)=H(x,a)$ (see
\cite{Wi96}). However, since $F_a$ has multiple factors, we get that
$D_y(F_a)=0$, so $H(x,a)=0$. Therefore, $z-a$ is a factor of
$D_y(F)$, and consequently (ii) occurs. Hence, in this case, the
result is proved. Now, let us see the case when $R=0$. Since
$\deg_y(F)>0$, then by Theorem \ref{main-theorem} it holds that
$M\neq 0$ but $M\in {\Bbb R}[z]$. Hence, by Theorem
\ref{main-theorem}, ${\mathcal A}$ is the set of real roots of $M$,
so clearly if $a$ satisfies (ii), then $a\in {\mathcal A}$. Thus,
let us see (i). For this purpose, let $A(z)$ be the leading
coefficient of $F$ w.r.t. $y$, and let $a$ be a real root of $A$. If
$\deg_y(F)>1$, then reasoning like before we get that $z-a$ is a
factor of $M$, and hence $a\in {\mathcal A}$. If $\deg_y(F)=1$, then
$\frac{\partial F}{\partial y}=A(z)$, and $\Res_y(F,\frac{\partial
F}{\partial y})=A(z)$. Hence, in this case it also holds that $z-a$
is a factor of $\Res_y(F,\frac{\partial F}{\partial y})$, and
consequently of $M$. Therefore, we also get that $a\in {\mathcal
A}$. In order to prove (iii), we proceed as in the case $R\neq 0$.
\qed

Now, we can prove the following result concerning to the level
curves of $S$.

\begin{theorem} \label{charact-reality}
Let ${\mathcal A}$ be a critical set of the surface $S$ determined
by applying Theorem \ref{main-theorem}. Then, $S$ is real if and
only if there exists at least one real level curve $S_a$ of $S$,
with $a\in {\Bbb R}$ and $a\notin {\mathcal A}$.
\end{theorem}

{\bf Proof: }If $S$ is real, then, by Theorem \ref{th-Lang}, there
exists a regular real point $P\in S$. Thus, the implication
$(\Rightarrow)$ follows from Implicit Function Theorem. Let us
consider the implication $(\Leftarrow)$. For this purpose, we
separately analyze the cases when $F \in {\Bbb R}[x,z]$, and when
$F$ depends on the variable $y$. We start with the case $F\in{\Bbb
R}[x,z]$. In this situation, let ${\mathcal C}_{xz}$ be the plane
algebraic curve defined by $F$ in the $xz$-plane. Thus, $S$ is real
iff ${\mathcal C}_{xz}$ is real. So, in order to prove that $S$ is
real, it suffices to prove that ${\mathcal C}_{xz}$ is. Let us see
that this holds. Now, by hypothesis there exists $a\in {\Bbb R}$
such that $a\notin {\mathcal A}$, and verifying that the
corresponding level curve $S_a$ is real. Since $F$ does not depend
on $y$, the level curves of $S$ are lines normal to the $xz$-plane.
Therefore, since $S_a$ is a real curve, one has that the
intersection point of $S_a$ with the $xz$-plane, which we denote as
$P_a$, is also real. By Theorem \ref{main-theorem}, ${\mathcal A}$
is the set of real roots of the discriminant $D_x(F)$. Thus, $a$ is
not a root of $D_x(F)$. Therefore, $P_a$ is not a singular point of
${\mathcal C}_{xz}$ and consequently ${\mathcal C}_{xz}$ is a real
curve. Thus, $(\Leftarrow)$ holds for the case when $F=F(x,z)$.
Finally, let us see that $(\Leftarrow)$ also holds when $F$ depends
on the variable $y$, i.e. when $\deg_y(F)>0$. In order to see this,
let $S_a$ be a level curve of $S$, real, and corresponding to the
intersection of $S$ with the real plane $z=a$, where $a\notin
{\mathcal A}$. Since $a\notin {\mathcal A}$, by Proposition
\ref{prop-1} the polynomial $F_a(x,y)=F(x,y,a)$ is square-free.
Thus, since $S_a$ is real and $F_a(x,y)$ is square-free, we have
that $S_a$ has at least one real non-singular point $(x_a,y_a)\in
S_a$. Then, $(x_a,y_a,a)$ is a real non-singular point of the
surface $S$ (note that $\nabla(F_a)=(F_x(x,y,a),F_y(x,y,a))$), so by
Theorem \ref{th-Lang} the surface $S$ is real. \qed

This theorem can be used to derive an algorithm for checking the
reality of an algebraic surface. For this purpose, note that the
condition in Theorem \ref{charact-reality}, i.e. the existence of a
real level curve of $S$ corresponding to a non-critical level value,
can be tested by checking the reality of the level curves
corresponding to intermediate $z$-values in between two consecutive
critical level values. More precisely, one has the following
algorithm:

{\sf \underline{Algorithm:} (Reality of an algebraic surface $S$)}
{\sf \underline{Given}} an algebraic surface $S$ implicitly defined
by a real polynomial $F(x,y,z)$, square-free, with no factor only
depending on the variable $z$, and such that $\lcoeff_y(F)$ does not
depend on the variable $x$, the algorithm {\sf \underline{decides}}
whether $S$ is real.

\begin{itemize}
\item [(1)] Compute a critical set of $S$, ${\mathcal A}=\{a_1,\ldots,a_r\}$
by means of Theorem \ref{main-theorem}. Let $a_0=-\infty$,
$a_{r+1}=\infty$.
\item [(2)]Check whether there exists $i\in
\{0,\ldots,r\}$ such that the plane algebraic curve defined by
$F(x,y,\xi_i)$, where  $\xi_i$ is taken in the interval
$(a_i,a_{i+1})$, is real. If it is, then {\bf return} $\ll S$ is
real$\gg$ {\bf else  return} $\ll S$ is not real$\gg$.
\end{itemize}

\begin{remark}
Note also that in case that the surface is not real there would be
three alternatives: (i) the real part of the surface reduces to a
space curve; (ii) it consists of finitely many points;  (iii) it is
empty. Moreover, one can algorithmically decide which is the case by
inspecting the level curves. More precisely, in case (iii) all the
level curves are empty; in case (ii), there are just finitely many
non-empty level curves, all of them corresponding to $z$-critical
level values, and consisting of finitely many real points. Finally,
case (i) is identified when (ii) and (iii) do not happen, and all
the level curves corresponding to non-critical $z$-values are either
empty or consisting in finitely many real points.
\end{remark}

\begin{example}\label{S1}
Let $S$ be the algebraic surface defined by \[F(x,y,z)=
(x^2-1)^2+(y^2-1)^2+(z^2-1)^2-3/2.\] Note that $S$ satisfies the
imposed hypotheses. Let us see whether $S$ is real. For this
purpose, we apply Theorem \ref{main-theorem} to obtain the following
$z$-critical level set of $S$:
\[\begin{array}{c}{\mathcal A}_z=\{-1.491557867, -1.306562965,
-0.5411961001,\\
0.5411961001, 1.306562965, 1.491557867\}\end{array}\] Now, we check
if there exists some level curve, corresponding to a $z$-value not
in ${\mathcal A}_z$, which is real. For $z<-1.491557867$ we get that
the level curves are empty over ${\Bbb R}^2$, but for $z=-7/5$,
which is intermediate between $-1.491557867$ and $-1.306562965$, we
get the $z$-slice
\[
\begin{array}{c}
\{(x^2-1)^2+(y^2-1)^2-723/1250=0,z=-7/5\}
\end{array}
\]
which is real. Therefore, we conclude that $S$ is real (see Figure
1).
\end{example}

\begin{example}\label{S2}
Consider the surface $S$ defined by \[F(x,y,z)=
x^4+2x^2y^2-2x^2+y^4-2y^2+1+z^2.\] Note that $S$ satisfies the
imposed hypotheses. In this case, a $z$-critical set of $S$ is
${\mathcal A}_z=\{0\}$, i.e. the $z$-slices of $S$ have at most
three different topology types. However, the $z$-slices for $z=-1$
and $z=1$ are empty curves over ${\Bbb R}$, so from Theorem
\ref{main-theorem} we deduce that for $z>0$ and $z<0$ the surface is
empty over the reals. Therefore, $S$ is not real. In fact, the only
real points of the surface are the points of the $z$-slice
corresponding to $z=0$, which is the circle
\[
\begin{array}{c}
\{(x^2+y^2-1)^2=0,z=0\}
\end{array}
\]
\end{example}

\section{Second application: Compactness}\label{compactness}

Here we show how to use level curves to algorithmically decide
whether $S$ is compact. Now, since $S$ is implicitly defined by a
polynomial $F\in {\Bbb R}[x,y,z]$, then it is obviously closed.
Thus, in order to check whether it is compact, it suffices to check
whether it is bounded, which is equivalent to decide whether it is
bounded w.r.t. the $x$, $y$ and $z$ variables, respectively.

Let us see how to check whether $S$ is bounded w.r.t. the variable
$z$. For this purpose, let ${\mathcal A}=\{a_1,\ldots,a_r\}$ be a
critical level set of $S$, where $a_1<\cdots<a_r$. Then, $S$ is
bounded w.r.t. the variable $z$ iff for $z>a_r$ and $z<a_1$, the
level curves of $S$ are empty over ${\Bbb R}^2$. Moreover, since by
Theorem \ref{main-theorem} the topology type of the level curves of
$S$ stays invariant for $z>a_r$, and also for $z<a_1$, in order to
check whether the condition holds it suffices to take $z_0<a_1$ and
$z_{r+1}>a_r$, and then to analyze whether the level curves
$S_{z_0},S_{z_{r+1}}$ are empty or not over ${\Bbb R}^2$. Note that
for this purpose one may adapt the strategy for deciding whether a
given algebraic curve is real.

Similarly for the $x$ and $y$ variables. However, observe that in
order to compute $\xi$-critical sets, with $\xi\in \{x,y,z\}$, by
means of Theorem \ref{main-theorem}, one needs that the hypotheses
of Theorem \ref{main-theorem} hold not only for the variable $z$,
but also for $x$, $y$, respectively. To ensure that this happens,
one may always apply if necessary a linear transformation so that
the polynomial $F\in {\Bbb R}[x,y,z]$ defining $S$ has no univariate
factor, and verifies that $\lcoeff_x(F)$, $\lcoeff_y(F)$ and
$\lcoeff_z(F)$ are all constant. Observe that this kind of
transformations preserves the topological properties of the surface.

Thus, one may derive the following algorithm:

{\sf \underline{Algorithm:} (Compactness of an algebraic surface
$S$)} {\sf \underline{Given}} an algebraic surface $S$ implicitly
defined by a real polynomial $F(x,y,z)$, square-free, with no
univariate factor, and such that
$\lcoeff_x(F),\lcoeff_y(F),\lcoeff_z(F)$ are constant, the algorithm
{\sf \underline{decides}} whether $S$ is compact.

\begin{itemize}
\item [(1)] Compute an upper bound $k_z$ of the absolute value of
the elements of a critical set of $S$.
\item [(2)] Check if the plane algebraic curves defined by
$F(x,y,-k_z-1)$ and $F(x,y,k_z+1)$ are both empty over ${\Bbb R}^2$.
If this does not happen, then {\bf return} $\ll S$ is not
compact$\gg$.
\item [(3)] Proceed in an analogous way for the variables $x$ and
$y$. If all the tested plane curves  are empty over ${\Bbb R}^2$,
{\bf return} $\ll S$ is compact$\gg$,  {\bf else return} $\ll S$ is
not compact$\gg$.
\end{itemize}

\begin{remark} \label{bounds}
In this case one does not need to compute  the real roots of the
polynomial provided by Theorem \ref{main-theorem}, but just upper
and lower bounds on them. This can be done by applying existing
algorithms (see for example \cite{Mignotte}).
\end{remark}

\begin{example}\label{S3}
Let $S$ be the surface in Example \ref{S1}, which fulfills all the
requirements of the algorithm before, and let us see whether it is
compact. For this purpose, we have to check whether it is bounded.
This is equivalent to checking whether the $z$-slices below the
least $z$-critical level value and above the greatest $z$-critical
level value are both empty over ${\Bbb R}^2$; similarly for $y$ and
$x$. In this sense, for $z=-2$ and $z=2$ one gets
\[
\begin{array}{c}
\{(x^2-1)^2+(y^2-1)^2+15/2=0,z=-2\}\\
\{(x^2-1)^2+(y^2-1)^2+15/2=0,z=2\}
\end{array}
\]
which are obviously empty over ${\Bbb R}^2$, so the condition holds
for $z$. By symmetry it also holds for $y$ and $x$. Therefore, we
deduce that $S$ is bounded, and therefore it is compact (see Figure
1).
\end{example}

\section{Third application: Plotting Boxes}\label{plotting}

Here we address the problem of computing an interval $I=[-a,a]\times
[-b,b]\times [-c,c]\subset {\Bbb R}^3$, so that the plotting of $S$
in $I$ shows the main relevant topological features of $S$. For this
purpose, the information on the $\xi$-level curves of $S$, $\xi\in
\{x,y,z\}$, is used.

More precisely, we consider the following definition, which provides
a criterion to compute $I$.

\begin{definition}
We say that the interval $[-m_x,m_x]\times [-m_y,m_y] \times
[-m_z,m_z]\subset {\Bbb R}^3$ is {\sf suitable for plotting $S$} if,
for $\xi\in \{x,y,z\}$,  $-m_{\xi},m_{\xi}$ are not $\xi$-critical
level values of $S$ and $[-m_{\xi},m_{\xi}]$ contains {\sf all} the
$\xi$-critical level values of $S$.
\end{definition}

Thus, if $I$ is ``suitable for plotting" $S$, one can be sure that
out of $I$ there is no change in the topology type of the
$\xi$-level curves of $S$. Note that the computation of a suitable
$I$ requires to compute critical level sets for the variables
$x,y,z$, respectively, so one requests the same hypotheses as in
Section \ref{compactness}. Observe also that Remark \ref{bounds}
also holds for this case. Thus, the following algorithm is derived:

{\sf \underline{Algorithm:} (Suitable interval for plotting an
algebraic surface $S$)} {\sf \underline{Given}} an algebraic surface
$S$ implicitly defined by a real polynomial $F(x,y,z)$, square-free,
with no univariate factor, and such that $\lcoeff_x(F)$,
$\lcoeff_y(F)$, $\lcoeff_z(F)$ are constant, the algorithm {\sf
\underline{determines}} a suitable interval $I\subset {\Bbb R}^3$
for plotting $S$.
\begin{itemize}
\item[(1)] For $\xi\in\{x,y,z\}$ compute an upper bound $k_\xi$ of the absolute
values of the elements of a $\xi$-critical set.
\item [(2)] {\bf Return} the interval $I=[-k_x-1,k_x+1]\times
[-k_y-1,k_y+1]\times [-k_z-1,k_z+1]$.
\end{itemize}

\begin{example}\label{S4}
Consider  again the surface in Example \ref{S1}. Here, we have that
\[\begin{array}{c}{\mathcal A}_z=\{-1.491557867, -1.306562965,
-0.5411961001,\\
0.5411961001, 1.306562965, 1.491557867\}\end{array}\]is a
$z$-critical level set of $S$. Furthermore, by symmetry,\[{\mathcal
A}_z={\mathcal A}_x={\mathcal A}_y\]Thus, the interval
\[I=[-1.5,1.5]\times[-1.5,1.5]\times [-1.5,1.5]\]is suitable for
plotting $S$. The picture of the part of $S$ lying in $I$ is shown
in Figure 1.
\end{example}

\begin{figure}[ht]
\begin{center}
\centerline{
\psfig{figure=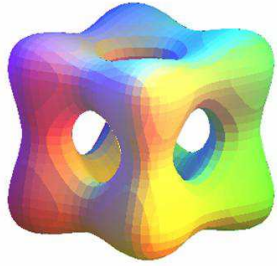,width=8cm,height=6cm}}
\end{center}
\caption{$F(x,y,z)=(x^2-1)^2+(y^2-1)^2+(z^2-1)^2-3/2$}
\end{figure}

\section{Fourth application: Connectivity}\label{plotting}

There exist purely symbolic algorithms, based for instance on the
notion of {\sf roadmap}, to compute the connected components of
$S$ (see \cite{Basu} for further details. Here we propose an
alternative approach to solve this question by means of a
symbolic-numeric algorithm based on level curves. The empirical
performance of this new method is quite satisfactory. The main
idea of the algorithm is to determine how the connected components
of a level curve corresponding to a non-critical $z$-value join to
the connected components of the level curves corresponding to the
critical level values immediately below and above, respectively.
For this purpose, we take points on each connected component of
the non-critical $z$-value, and then we generate a space curve (as
the solution of a system of differential equations) which connects
it with some connected component of the critical $z$-level curve
immediately below/above (see Figure 2). Thus, once we know how to
join the connected components of the level curves corresponding to
non-critical and critical $z$-values, the number of connected
components of $S$ can be obtained as the number of ``connected
chains" (whose elements are connected components of level curves)
computed in the process.

In order to solve this problem, we require some more conditions on
the surface $S$ to be analyzed. More precisely, the following
hypotheses must be satisfied:

\begin{itemize}
\item [(i)] $S$ is defined by an square-free polynomial $F$, with no
factor just depending on the variables $y,z$. If this holds, then
$\gcd(F,F_x)=1$, so the variety ${\mathcal C}$ defined by $F=F_x=0$
is a space algebraic curve; recall that $F_x$ denotes the partial
derivative of $F$ w.r.t. $x$.
\item [(ii)] There does not exist any plane $z-a=0$
containing infinitely many points of the curve ${\mathcal C}$
introduced in (i), i.e. not containing infinitely many points of $S$
where $F_x$ vanishes.
\item [(iii)] $S$ is not asymptotic to any plane of equation $z-b=0$, where
$b$ is a critical level value.
\end{itemize}

Observe that in case that $S$ does not fulfill some of these
conditions, almost all linear transformations lead to a new surface
where the three requirements hold. Note that linear transformations
preserve the topological features of the surface.

Now, in the sequel let ${\mathcal A}=\{a_1,\ldots,a_r\}$, where
$a_1<\cdots<a_r$, be a $z$-critical level set of $S$; furthermore,
we set $a_0=-\infty$, $a_{r+1}=+\infty$. Moreover, let
$b_0<\cdots<b_r$ verify $a_i<b_i<a_{i+1}$ for all $i\in
\{0,\ldots,r\}$. With this notation, our problem is to decide, for
each $a_i,b_i,a_{i+1}$, how to connect the connected components of
the level curves $S_{b_i}, S_{a_i}$, and $S_{b_i}, S_{a_{i+1}}$,
respectively. Observe that  computing the topology graph of a level
curve (which is planar) one obtains the connected components of it
(see for instance \cite{Lalo}). Also, note that it may happen that
the level curves for values in some of the intervals $(a_i,b_i)$ are
empty over the reals, in which case there is no need of connecting
$S_{b_i}, S_{a_i}$, and $S_{b_i}, S_{a_{i+1}}$. In addition, observe
that, since $S$ verifies condition (iii), every connected component
of some level curve of $S$ corresponding to a non-critical
$z$-value, joins to some connected component of the level curve
corresponding to the $z$-critical level value immediately below
(resp. above). In fact, this is the reason why we request condition
(iii).

For this purpose, the strategy is to use a symbolic-numeric
algorithm which essentially works as follows:

{\sc \underline{Algorithm:} (Connected components of an algebraic
surface $S$)} {\sf \underline{Given}} an algebraic surface $S$
implicitly defined by a real polynomial $F(x,y,z)$ verifying the
conditions (i), (ii), (iii), the algorithm {\sf
\underline{determines}} the number of connected components of $S$
and a description of them in terms of level curves.

\begin{itemize}
\item [(0)] Compute the topology graph of
$S_{b_i}$ (if $S_{b_i}\cap {\Bbb R}^{3}=\emptyset$ take another
$i$), and the singular points of $S_{a_i}$; the information on the
singular points of $S_{a_i}$ will be used at step (2), in some cases
(see Remark \ref{singu}), and also at step (3)).
\item [(1)] Take a real point $P$ in each connected component of
$S_{b_i}$.
\item [(2)] Use a path continuation method to connect $P$ with some
point $Q$ in $S_{a_i}$, to be computed by the algorithm; in order to
do this, we travel from $P$ to $Q$ by following a space curve,
contained in the surface, which is computed as the solution of a
system of differential equations.
\item [(3)] Identify the connected
component of $S_{a_i}$ where the final point $Q$ belongs to. For
this purpose, we compute the topology graph of $S_{a_i}$ introducing
the point $Q$ as a vertex of the graph.
\item [(4)] Join by an edge the starting connected component of $S_{b_i}$ and the
finally reached connected component of $S_{a_i}$.
\item [(5)] Proceed in an analogous way to connect $P$ with some
point $Q^{\star}$ in $S_{a_{i+1}}$.
\item [(6)] After carrying out the computation for all the connected
components of all the $S_{b_i}$, one gets several connected chains,
each one consisting of some connected components of the $S_{b_i}$,
$S_{a_i}$ joined by edges. The number of connected chains, is the
number of connected components of the surface.
\end{itemize}

Observe that in the execution of step (0), points on each connected
component of $S_{b_i}$  are computed (see e.g. \cite{Lalo} or
\cite{Hong}), so step (1) can be executed afterwards. Now, let us
describe with more detail step (2). We consider the solution to the
following system of differential equations:
\[
\left \{\begin{array}{l}
x'=-F_y+\displaystyle{\frac{F_z}{F_x}}\\
y'=F_x\\
z'=-1 \\
x(0)=x_i;y(0)=y_i;z(0)=b_i
\end{array}\right.
\]
where $P:=(x_i,y_i,b_i)\in S_{b_i}$. In case that this differential
system has a symbolic solution $(x(t),y(t),z(t))$, one may see that
it corresponds to a space curve contained in $S$. Moreover, since
$z'(t)=-1$ and $S$ fulfills condition (iii), then this space curve
reaches $z=a_i$. Furthermore, since one may decompose the part of
$S$ with $z\in (a_i,b_i)$ into non-intersecting pieces, each one
corresponding to a different connected component of $S_{b_i}$ (see
\cite{JGRS} for a careful proof of this fact), the choice of the
initial point for a particular connected component of $S_{b_i}$ does
not affect the connected component finally reached. In other words,
the connected component reached at $z=a_i$ is always the same for
all the points of a same connected component of $S_{b_i}$ (see also
Figure 2).

\para

\para

\begin{figure}[ht]
\begin{center}
\centerline{ \psfig{figure=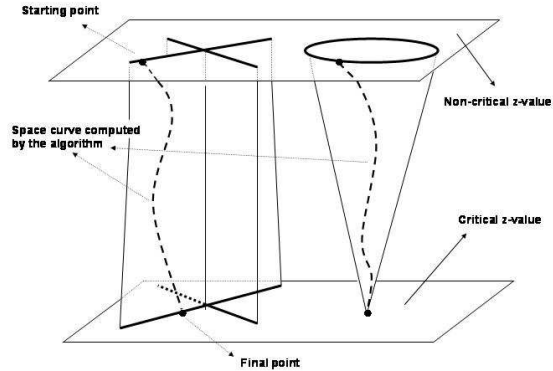,width=8cm,height=6cm}}
\end{center}
\caption{Idea of the connectivity algorithm}
\end{figure}

However, in general the differential system above may not have a
symbolic solution, so numerical methods must be applied. In our
case, we used the package of {\sc maple} for numerically integrating
differential equations. Here, one may see that as the numerical
method goes on, the error considerably grows, so the point of
$z=a_i$ finally reached cannot be recognized as belonging to any
connected component of $S_{a_i}$. For this purpose, at each step the
solution provided by numerical integration must be corrected. In
order to do this, each solution
$\tilde{P}_{i,k}=(\tilde{x}_{i,k},\tilde{y}_{i,k},\tilde{z}_{i,k})$
is corrected to $\bar{P}_{i,k}=(\bar{x}_{i,k},
\bar{y}_{i,k},\tilde{z}_{i,k})$ by computing a point of the level
curve $S_{\tilde{z}_{i,k}}$ close to the solution
$(\tilde{x}_{i,k},\tilde{y}_{i,k},\tilde{z}_{i,k})$ (observe that
the $z$-coordinate is the same in $\tilde{P}_{i,k}$ and in
$\bar{P}_{i,k}$). For this purpose, we take the line passing through
$\tilde{P}_{i,k}$ in the direction of
$\vec{v}=(F_x(\tilde{P}_{i,k}),F_y(\tilde{P}_{i,k}))$, and we
compute the intersection points of this line with
$S_{\tilde{z}_{i,k}}$. This new point is used to go on with the
numerical integration process. Although we have not analyzed the
convergence of the method, the empirical results we get in this way
are satisfactory. In this sense, the following remark must be taken
into account.

\begin{remark} \label{singu} If $K$ is a singular
point of $S_{a_i}$ (notice that these points are computed in the
initial step of the algorithm), then we assume that it is reached
whenever the distance between $\bar{P}_{i,k}$ and $K$ is smaller
than a sufficiently small $\epsilon$ previously fixed. Thus, in this
case the computation stops and we assume that the starting point of
$S_{b_i}$ is connected with $K$, i.e. that $Q=K$.
\end{remark}

In addition, there are two more situations which must be examined
carefully:
\begin{itemize}
\item If, before reaching the level plane $z=a_i$, the numerical
integration process hits a point where $F_x$ vanishes, then the
method fails. This situation can be prevented by detecting whether
$|F_x(\tilde{x}_{i,k},\tilde{y}_{i,k},\tilde{z}_{i,k})|<\epsilon$;
if this happens, we choose a different point of
$S_{\tilde{z}_{i,k}}$ to go on with the numerical process.
\item It may happen that $S$ contains a 1-dimensional subset
${\mathcal L}$ of singular points, where ${\mathcal L}$ is
``isolated" in the following sense: given any point $P\in {\mathcal
L}$, there exists a Euclidean neighborhood $E_p$ of $P$ such that
every point of $S\cap E_p\cap {\Bbb R}^{3}$ is also a point of
${\mathcal L}$. For example, the handle of the Whitney Umbrella
$x^2-y^2z=0$, which is obtained for negative values of $z$, provides
an example of this situation. Unless ${\mathcal L}$ is parallel to
the $xy$-plane, this phenomenon can be detected by identifying the
presence of isolated points in non-critical level curves.

Hence, assume that $S_{b_i}$ has some isolated point $B$, therefore
belonging to some ${\mathcal L}$. The points of ${\mathcal L}$ are
singular points of $S$, so $F_x$ vanishes at $B$ and therefore the
system of differential equations before cannot be used to compute
the point $A$ in $S_{a_i}$ which must be connected with $B$. Now, in
this case we use the projections of the curve ${\mathcal L}$ onto
two coordinates planes to compute how the points of ${\mathcal L}$
in $z=b_i$ and $z=a_i$, respectively, are connected, in analogy with
the method described in \cite{JG-Sendra} to compute the topology of
space algebraic curves. More precisely, we consider a coordinate
plane so ${\mathcal L}$ is not normal to it. We project ${\mathcal
L}$ onto this plane, we determine the projection $\pi(B)$ of the
point $B$, and then we determine the projection $\pi(A)$ of the
point $A$ by either applying a path continuation method (over the
projection), or by computing the part of the graph associated with
the projection between $z=b_i$ and $z=a_i$. Carrying out this
process on two coordinate planes, the point $A$ is obtained.




\end{itemize}

Finally, in order to connect $z=b_i$ and $z=a_{i+1}$, we apply an
analogous process to the differential equation system:
\[
\left\{\begin{array}{l}
x'=-F_y-\displaystyle{\frac{F_z}{F_x}}\\
y'=F_x\\
z'=1\\
x(0)=x_i;y(0)=y_i;z(0)=b_i
\end{array}\right.
\]
Note here that the third equation is different from the system
before ($z'=1$ instead of $z'=-1$), since in this case one has to
move ``up" from $z=b_i$ to $z=a_{i+1}$. One may see that also the
first equation has changed. However, the space curve that one
obtains by integrating these equations lies also in the surface $S$.

\begin{example}\label{S5}
Consider the algebraic surface $S$ defined by
$F(x,y,z)=x^2+y^2+z^2+2xyz-1$. A $z$-critical level set of $S$ is
${\mathcal A}_z=\{-1,1\}$. Because of the symmetry of the surface,
we have that ${\mathcal A}_x={\mathcal A}_y={\mathcal A}_z$, so
for example $[-2,2]\times[-2,2]\times[-2,2]$ is suitable for
plotting $S$. A plotting of $S$ in this interval can be seen in
Figure 3; this figure was computed with {\sc maple}.

\begin{figure}[ht]
\begin{center}
\centerline{
\psfig{figure=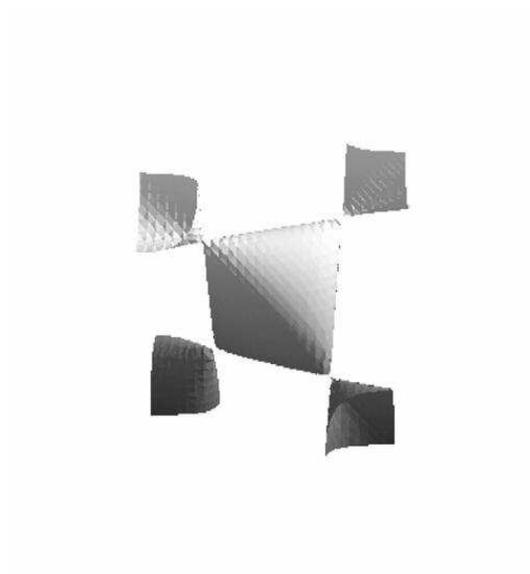,width=10cm,height=8cm}}
\end{center}
\caption{The cubic $x^2+y^2+z^2+2xyz-1=0$}
\end{figure}

Observe that, from Figure 3, it is not completely clear whether
the surface is connected or not. However, using the algorithm  we
check that $S$ is connected. Indeed, Figure 4 shows the different
topology types corresponding to the $z$-level curves for the cases
$z<-1,z=-1,-1<z<1,z=1,z>1$, respectively. Moreover, also in Figure
4, each connected component of each $z$-level curve has been
joined with the corresponding connected component of the
$z$-critical level curve immediately above/below, according to the
algorithm in this section. The connections between connected
components of level curves are represented by dotted lines. Here,
one may see that there is just one connected ``chain", formed by
the connected components of the level curves which are joined one
another. Hence, $S$ has just one connected component, and
therefore it is connected.

\para

\para

\begin{figure}[ht]
\begin{center}
\centerline{
\psfig{figure=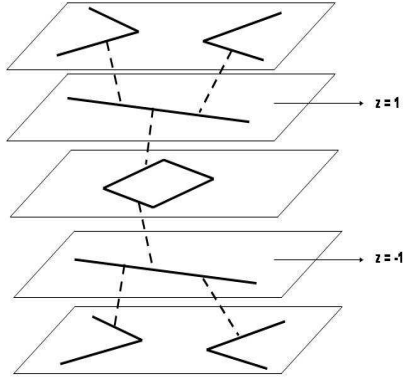,width=6cm,height=6cm}}
\end{center}
\caption{Connections between the level curves of
$x^2+y^2+z^2+2xyz-1=0$}
\end{figure}
\end{example}

\para

\begin{example}\label{S6}
Consider the quartic $(x^2-1)^2+(y^2-1)^2+(z^2-1)^2-3/4=0$. In this
case, a $z$-critical level set of $S$ is
\[{\mathcal
A}_z=\{-1/2-\sqrt{3}/2,1/2-\sqrt{3}/2,-1/2+\sqrt{3}/2,1/2+\sqrt{3}/2\}\]The
different topology types for the $z$-level curves are shown in
Figure 5; moreover, one may check that for $z>1/2+\sqrt{3}/2$ and
$z<-1/2-\sqrt{3}/2$, the level curves are empty over ${\Bbb R}^2$.
Furthermore, also in this picture we have represented (in dotted
lines) how the connected components of these level curves are
joined with each other. This information has been computed by
using the algorithm in this section.

\newpage

\begin{figure}[ht]
\begin{center}
\centerline{
\psfig{figure=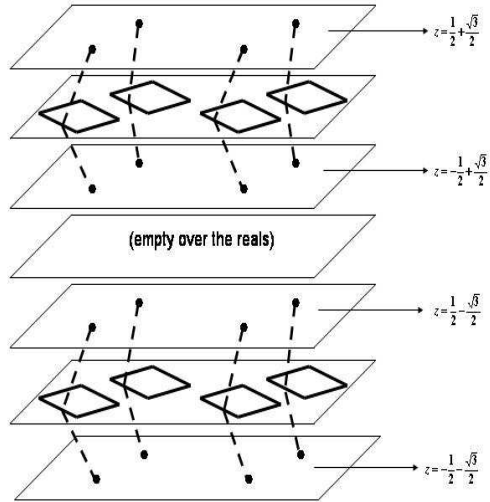,width=8cm,height=8cm}}
\end{center}
\caption{Connections between the level curves of
$(x^2-1)^2+(y^2-1)^2+(z^2-1)^2-3/4=0$}
\end{figure}

\para

Since in Figure 5 one may see 8 different connected chains of
components of level curves, one has that $S$ has 8 connected
components. A plotting of this quartic can be seen in Figure 6.

\begin{figure}[ht]
\begin{center}
\centerline{ \psfig{figure=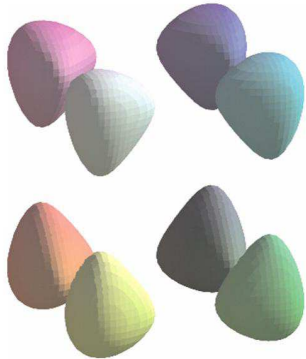,width=6cm,height=5cm}}
\end{center}
\caption{The quartic $(x^2-1)^2+(y^2-1)^2+(z^2-1)^2-3/4=0$}
\end{figure}
\end{example}

\end{document}